\documentclass[10pt]{article}
\usepackage{amssymb, amsmath, amsfonts,latexsym, oldgerm,amscd,amsthm}
\usepackage{relsize}
\usepackage{url}

\textwidth13cm \usepackage[all]{xy} 
\newcommand{\tn}{\textnormal} \newcommand{\lra}{\longrightarrow}
 \newtheorem{de}{Definition}[section]

 \newtheorem{pr}[de]{Proposition}
\newtheorem{thm}[de]{Theorem} \newtheorem{lem}[de]{Lemma}
 \newtheorem{nt}[de]{Notation}
 
 \newtheorem{rmk}[de]{Remark}

\def\E{{\rm E}}

\def\O{{\rm O}}
\def\EO{{\rm EO}}
\def\GL{{\rm GL}}

\def\Um{{\rm Um}}

\def\End{{\rm End}}
\def\Aut{{\rm Aut}}
\def\Trans{{\rm Trans}}
\def\ETrans{{\rm ETrans}}

\def\Hom{{\rm Hom}}

\def\lra{\longrightarrow}

\newcommand{\gm}{\mathfrak{m}}

\begin{document}

\title{Equality of orthogonal transvection group and elementary orthogonal transvection group}

\author{Pratyusha Chattopadhyay   \\
{\small Statistics and Mathematics Unit, Indian Statistical Institute, Bangalore 560 059, India}
}

\date{}

\maketitle

\begin{center}{\small {\it 2010 MSC: {11C20, 13C10, 14L35, 15A24, 15A54, 15A63, 20H25}}}
\end{center}

\begin{center}
\end{center}

{\small ~~~~Abstract: H. Bass defined orthogonal transvection group of an orthogonal module
and elementary orthogonal transvection group of an orthogonal module with a hyperbolic direct summand. 
We also have the notion of relative orthogonal transvection group and relative elementary orthogonal 
transvection group with respect to an ideal of the ring. According to the definition of Bass relative elementary 
orthogonal transvection group is a subgroup of relative orthogonal transvection group of an orthogonal 
module with hyperbolic direct summand. Here we show that these two groups are the same
in the case when the orthogonal module splits locally.}

\section{\large Introduction}

In Section 5 of \cite{SV} L.N. Vaserstein proved that first row of an elementary linear matrix of even size 
(bigger than or equal to 4) is the same as the first row of a symplectic matrix of the same size w.r.t. an alternating 
form. This result motivated us to prove that the orbit of a unimodular row of even size under the action of 
elementary linear group is same as the orbit of a unimodular row of same size under the action of elementary 
symplectic group (see Theorem 4.1, \cite{cr}); we also proved a relative (to an ideal of the ring) version of this result (see Theorem 5.5, \cite{cr}). Generalising this result in the setting of finitely generated projective modules
involving the transvection groups as defined by H. Bass,  we proved that 
in the case of a symplectic module with a hyperbolic direct summand the orbits of any unimodular element from the 
symplectic module, under the actions of elementary linear transvection group and the  elementary symplectic transvection group
coincide(see Theorem 6.1, \cite{cr2}).
While proving the above result on equality of orbits of unimodular 
elements of symplectic modules, we observed that in the relative case to an ideal of the ring, the equality holds between the linear transvection 
group and the elementary linear transvection group (see Proposition 4.10, \cite{cr2}).
We also noticed that in the
relative case to an ideal of the ring, the symplectic transvection group and the elementary symplectic transvection 
group coincide (see Theorem 5.23, \cite{cr2}). In the absolute case 
the  equalities for linear transvection group, symplectic
transvection group, and orthogonal transvection group  with the corresponding elementary transvection groups
were proved in \cite{bbr}. In view of the above results it  is  natural to ask whether the equality of orthogonal transvection group and elementary orthogonal 
transvection group holds in the relative case to an ideal of the ring. 
In this article we prove the equality of these two groups in the case when the orthogonal module splits locally.

\section{\large Preliminaries}

\medskip
In this article we will always assume that $R$ is a commutative ring with unit. 
A row ${ v} = (v_{1}, \ldots, v_{n}) \in R^{n}$ is said to be {\it
  unimodular} if there are elements $w_{1}, \ldots, w_{n}$ in $R$ such
that $v_{1}w_{1} + \cdots + v_{n}w_{n} = 1$. Um$_{n}(R)$ will denote
the set of all unimodular rows ${ v} \in R^{n}$. Let $I$ be an ideal
in $R$. We denote by ${\rm Um}_n(R, I)$ the set of all unimodular rows
of length $n$ which are congruent to $e_1 = (1, 0, \ldots, 0)$ modulo
$I$. (If $I = R$, then ${\rm Um}_n(R, I)$ is ${\rm Um}_n(R)$).

\begin{de} {\rm Let $P$ be a finitely generated projective $R$-module. 
An element $p \in P$ is said to be {\it unimodular} if there exists a 
$R$-linear map $\phi: P \to R$ such that $\phi(p)=1$. The collection 
of unimodular elements of $P$ is denoted by $\Um(P)$. 

Let $P$ be of the form $R \oplus Q$ and have an element of the form
$(1,0)$ which correspond to the unimodular element. An element $(a,q)
\in P$ is said to be {\it relative unimodular} w.r.t. an ideal $I$ of 
$R$ if $(a,q)$ is unimodular and $(a,q)$ is congruent to $(1,0)$
modulo $IP$.  The collection of all relative unimodular elements
w.r.t. an ideal $I$ is denoted by ${\rm Um}(P,IP)$.  }
\end{de}

Let us recall that if $M$ is a finitely presented $R$-module and $S$
is a multiplicative set of $R$, then $S^{-1} {\rm Hom}_R(M,R) \cong
{\rm Hom}_{R_S}(M_S, R_S)$ (Theorem 2.13", Chapter I, \cite{lam}). 
Also recall that if $f=(f_1, \ldots,
f_n)\in R^n := M$, then $\Theta_M(f)=\{ \phi(f): \phi \in {\rm
  Hom}(M,R) \}= \sum_{i=1}^n Rf_i$. Therefore, if $P$ is a finitely
generated projective $R$-module of rank $n$, $\gm$ is a maximal ideal
of $R$ and $v\in \Um(P)$, then $v_\gm \in \Um_n(R_\gm)$. Similarly if
$v \in \Um(P,IP)$ then $v_\gm \in \Um_n(R_\gm, I_\gm)$.

\begin{de} {\bf Elementary Linear Group:} 
{\rm Elementary linear group E$_{n}(R)$ denote the subgroup of SL$_{n}(R)$ consisting of all
{\it elementary} matrices, i.e. those matrices which are a finite
product of the {\it elementary generators} E$_{ij}(\lambda) = I_{n} +
e_{ij}(\lambda)$, $1 \leq i \neq j \leq n$, $\lambda \in R$, where
$e_{ij}(\lambda) \in$ M$_{n}(R)$ has an entry $\lambda$ in its $(i,
j)$-th position and zeros elsewhere.}
\end{de}

In the sequel, if $\alpha$ denotes an $m \times n$ matrix, then we let
$\alpha^t$ denote its {\it transpose} matrix. This is of course an
$n\times m$ matrix. However, we will mostly be working with square
matrices, or rows and columns.

\begin{de} {\bf The Relative Groups E$_n(I)$, E$_n(R,I)$:} {\rm Let $I$ be
  an ideal of $R$. The relative elementary linear group {\rm E}$_n(I)$ is the
  subgroup of {\rm E}$_n(R)$ generated as a group by the elements {\rm
    E}$_{ij}(x)$, $x \in I$, $1 \leq i \neq j \leq n$.

  The relative elementary linear group ${\rm E}_n(R, I)$ is the normal
  closure of {\rm E}$_n(I)$ in {\rm E}$_n(R)$.

$($Equivalently, ${\rm E}_n(R, I)$ is generated as a group by ${\rm
    E}_{ij}(a) {\rm E}_{ji}(x)$E$_{ij}(-a)$,\tn{ with} $a \in R$, $x
  \in I$, $i \neq j$, provided $n \geq 3$ {\rm (see \cite{V3}, Lemma 8)}$)$.}
\end{de}

\begin{de} 
{\rm E$_n^1(R, I)$ is the subgroup of $\E_n(R)$ generated by the
  elements of the form $E_{1i}(a)$ and $E_{j1}(x)$, where $a \in R, x
  \in I$, and $2 \le i, j \le n$.  }
\end{de}

\begin{rmk} \label{unimodularOverLocalRing}
It is easy to check that if $v \in \Um_n(R,I)$, where $(R, \gm)$ is a
local ring and $I$ be an ideal of $R$, then $v = e_1 \beta$, for some
$\beta \in \E_n(R,I)$.
\end{rmk}

\begin{de} {\bf Orthogonal Group:} {\rm The orthogonal group $\O_{2n}(R)$ 
with respect to the standard symmetric matrix $\widetilde {\psi}_n = \underset{i=1}{\overset{n}\sum}
  e_{2i-1,2i} + \underset{i=1}{\overset{n}\sum} e_{2i,2i-1}$ is the collection
$\{\alpha \in {\rm GL}_{2n}(R)\,\,|\,\, \alpha^t \widetilde{\psi}_n
  \alpha = \widetilde{\psi}_n \}$. For an ideal $I$ of $R$, $\O_{2n}(R, I)$ 
  represents the kernel of the natural map $\O_{2n}(R) \lra \O_{2n}(R/I)$.} 
\end{de}

\tn{Let $\sigma$ denote the permutation of the natural numbers $\{1, 2, \ldots, 2n \}$ given
  by $\sigma(2i)=2i-1$ and $\sigma(2i-1)=2i$}.

\begin{de} \label{2.4} {\bf Elementary Orthogonal Group:} 
{\rm As in \S 2 of \cite{SK} we define for $z \in R$, $1\le i\ne j\le 2n$,
\begin{eqnarray*}
oe_{ij}(z) =  1_{2n} + z e_{ij} - z e_{\sigma (j) \sigma (i)} & {\rm if} ~ i\ne \sigma(j).
\end{eqnarray*}

It is easy to check that all these elements belong to O$_{2n}(R)$.
We call them {\it elementary orthogonal matrices} with respect to the standard
symmetric matrix $\widetilde{\psi}_n$ over $R$ and the
subgroup of {\rm O}$_{2n}(R)$ generated by them is called the
elementary orthogonal group {\rm EO}$_{2n}(R)$ with respect to the standard
symmetric matrix $\widetilde{\psi}_n$.}
\end{de}

\begin{de} {\bf The Relative Group EO$_{2n}(I)$, EO$_{2n}(R,I)$:}
{\rm Let $I$ be an ideal of $R$. The relative elementary group $\EO_{2n}(I)$ 
is the subgroup $\EO_{2n}(R)$ generated as a group by the elements
$oe_{i j}(x), x \in I$ and $1 \le i \ne j \le 2n$.

The relative elementary group $\EO_{2n}(R, I)$ is the normal closure of
$\EO_{2n}(I)$ in $\EO_{2n}(R)$.}
\end{de}

\begin{lem} \label{equiv-defn}
$\EO_{2n}(R, I)$ is  generated as a group by the elements of the form $g ~oe_{ij}(x) g^{-1}$, where
$g \in \EO_{2n}(R), x \in I$, and either $i=1$ or $j=1$.
\end{lem}

Proof: An element of the form $g ~oe_{ij}(x) g^{-1} \in \EO_{2n}(R, I)$, where
$g \in \EO_{2n}(R), x \in I$, and either $i=1$ or $j=1$. Consider an elementary generator 
$oe_{kl}(a) oe_{ij}(x) oe_{kl}(-a)$ of
$\EO_{2n}(R, I)$, where $a \in R, x \in I$, and $i, j \ne 1$. Then
\begin{eqnarray*}
&& oe_{kl}(a) oe_{ij}(x) oe_{kl}(-a)\\
&=& ^{oe_{kl}(a)} [oe_{i1}(x), oe_{1j}(1)] \\
&=& ^{oe_{kl}(a)} \{ oe_{i1}(x) ~  oe_{1j}(1) ~  oe_{i1}(-x) ~  oe_{1j}(-1) \} \\
&=& ^{oe_{kl}(a)} oe_{i1}(x)  ^{oe_{kl}(a)} \{ oe_{1j}(1) ~  oe_{i1}(-x) ~  oe_{1j}(-1) \}
\end{eqnarray*}
and hence the lemma follows.
\hfill{$\square$}

\begin{de}\label{defn-EO^1}
  {\rm The group $\EO_{2n}^1(R,I)$ is the subgroup of $\EO_{2n}(R)$
  generated by the elements of the form $oe_{1 i}(a)$ and $oe_{j
    1}(x)$, where $a \in R$, $x \in I$ and $3 \le i, j \le 2n$.}
\end{de}

In the following two lemmas we obtain some useful facts regarding elementary
orthogonal groups. An analogous result in the linear case was proved in 
\cite{vdK1} and in the symplectic case was proved in the Appendix
of \cite{cr2}.

\begin{lem} \label{subset}
Let $R$ be a commutative ring and $I$ be an ideal of $R$. Then $\EO_{2n}(R, I) 
\subseteq \EO_{2n}^1(R, I)$, for $n \ge 3$.
\end{lem}

Proof: Let us define $S_{ij} = \{ oe_{ij}(a) oe_{ji}(x) oe_{ij}(-a) : a \in R, x \in I \}$. It suffices 
to show that $\EO_{2n}^1(R, I)$ contains the set $S_{ij}$ for all $1 \le i \ne j \le 2n$. 
Note that $S_{ij} = S_{\sigma(j) \sigma(i)}$ and 
$S_{1j} \subseteq \EO_{2n}^1(R, I)$, for $3 \le j \le 2n$. First we state the following 
identities
\begin{align}\label{commprop}
[gh,k] &~ = ~ \big({}^g[h,k]\big)[g,k],\\
[g,hk] &~ = ~ [g,h]\big({}^h[g,k]\big),\\
{}^g[h,k] &~ = ~ [{}^gh,{}^gk],
\end{align}
where $^gh$ denotes $ghg^{-1}$ and $[g,h]=ghg^{-1}h^{-1}$. Using these
identities and the commutator law $[oe_{ik}(a), oe_{kj}(b)] = oe_{ij}(ab)$
we establish the inclusion. 

Note that $oe_{1j}(x), oe_{i1}(x) \in \EO_{2n}^{1}(R,I)$, for $3 \le
i,j \le 2n$ and $x \in I$. For $3 \le i,j \le 2n$ and $x \in I$, we have $oe_{ij}(x) = [oe_{i1}(x),
  oe_{1j}(1)] \in \EO_{2n}^{1}(R,I)$. In
      the following computation we will express the generators of
      $S_{ij}$ in terms of $oe_{i1}(x)$ or $oe_{1j}(a)$, where $x \in
      I, a \in R $. Also, note that we will use $\circledast$ to represent
      elements of $\EO_{2n}^1(R, I)$.

\begin{eqnarray*}
^{oe_{ij}(a)} oe_{ij}(x) &=& ^{oe_{ij}(a)} [oe_{j2}(1), oe_{2i}(x)] \\
&=& [oe_{i2}(a) oe_{j2}(1), oe_{2i}(x) oe_{2j}(ax)] \\
&=& ^{oe_{i2}(a)} [oe_{j2}(1), oe_{2i}(x) oe_{2j}(ax)] [oe_{i2}(a), oe_{2i}(x) oe_{2j}(ax)] \\
&=& ^{oe_{i2}(a)} [oe_{j2}(1), oe_{2i}(x)] ^{oe_{i2}(a) oe_{2i}(x)} [oe_{j2}(1), oe_{2j}(ax)] \\
&&  [oe_{i2}(a), oe_{2i}(x)] oe_{2j}(a^2x^2) oe_{2i}(x) \\
&=& ^{oe_{i2}(a)} oe_{ji}(x) ^{oe_{i2}(a)}[oe_{j2}(1) oe_{ji}(-x), oe_{2j}(ax)]  ~ \circledast\\
&=& \circledast ~ [oe_{j2}(1) oe_{ji}(x) oe_{j2}(ax), oe_{ij}(a^2x) oe_{2j}(ax)] ~ \circledast\\
&=& \circledast  ~ ^{oe_{j2}(1)} [oe_{ji}(x) oe_{j2}(ax), oe_{ij}(a^2x) oe_{2j}(ax)] \\
&& [oe_{j2}(1), oe_{ij}(a^2x) oe_{2j}(ax)] ~ \circledast \\
&=& \circledast ~ [oe_{ij}(x) oe_{j2}(ax), oe_{ij}(a^2x) oe_{i2}(-a^2x) ^{oe_{j2}(1)} oe_{2j}(ax)] \\
&& oe_{i2}(-a^2x) ^{oe_{ij}(a^2 x^2)} [oe_{j2}(1), oe_{2j}(ax)] ~ \circledast
\end{eqnarray*}

Since $S_{i2}, S_{j2} \subseteq
\EO_{2n}^1(R,I)$, therefore $S_{ij} \subseteq \EO_{2n}^1(R, I)$, for $3 \le i, j \le 2n$. Similarly $S_{ik}, 
S_{1k} \subseteq \EO_{2n}^1(R, I)$ will imply that $S_{i1} \subseteq \EO_{2n}^1(R, I)$, for $3 \le i \le 2n$.
\hfill{$\square$}

\begin{pr} \label{vanderk1} 
Let $R$ be a commutative ring and let $I$ be an
ideal of $R$. Then for $n \ge 3$ the following sequence is exact

\[ 1 \lra \EO_{2n}(R,I) \lra \EO^1_{2n}(R,I) \lra \EO^1_{2n}(R/I,0) \lra 1. \]

Thus $\EO_{2n}(R,I)$ equals $\EO^1_{2n}(R,I) \cap \O_{2n}(R,I)$.
\end{pr}

Proof: Let $f : \EO_{2n}^1(R, I) \lra \EO_{2n}^1(R/I, 0)$. Note that 
$\ker(f) = \EO_{2n}^1(R,I) \cap \O_{2n}(R,I)$. We shall prove that $\ker(f) = \EO_{2n}(R,I)$. 
Let $E = \prod_{k=1}^r
oe_{j_k 1}(x_k) oe_{1i_k}(a_k)$ be an element in the $\ker(f)$ and $E$ can be
written as $oe_{j_1 1}(x_1) \prod_{k=2}^r \gamma_k oe_{j_k 1}(x_k)
\gamma_k^{-1}$, where $\gamma_l$ is equal to
$\prod_{k=1}^{l-1}~oe_{1i_k}(a_k) \in \EO_{2n}(R)$, and hence $\ker(f)
\subseteq \EO_{2n}(R,I)$.  The reverse inclusion follows from the fact
that $\EO_{2n}(R,I) \subseteq \EO_{2n}^1(R,I)$ (see Lemma \ref{subset}). 
\hfill{$\square$}

\section{{\large Local Global Principle for Relative Elementary Group}}

In this section we prove Lemma \ref{ness4-E1} and Lemma \ref{E1-dil-strong} which will be used in proving the main result in the final section.
In Lemma \ref{E1-dil-strong} we obtain the Local-Global Principle for an extended
ideal for slightly larger group EO$_{2n}^1(R,I)$ than the relative group EO$_{2n}(R, I)$. This group was introduced in the linear case
by W. van der Kallen in \cite{vdK1}.
The Local-Global Principle for an extended
ideal in the linear, orthogonal and symplectic groups was proved in \cite{acr}.

The line of arguments given in the proofs below closely follow that in the  Local-Global principle for an extended ideal in the symplectic case 
in \cite{cr2}. However, as far as the the computational details are concerned, there are substantial deviations from \cite{cr2} in many steps.

\begin{lem} \label{ness4-E1} Let $R$ be a commutative ring and $I$ be an 
ideal of $R$. Let $n \ge 3$. Let $\varepsilon =
  \varepsilon_1 \ldots \varepsilon_r$ be an element of $\EO_{2n}^1(R,I)$,
  where each $\varepsilon_k$ is an elementary generator.  If
  $oe_{ij}(Xf(X))$ is an elementary generator of $\EO_{2n}^1(R[X],I[X])$,
  then
\begin{eqnarray*}
\varepsilon ~ oe_{ij}(Y^{4^r}Xf(Y^{4^r}X)) ~ \varepsilon^{-1} 
&=& \prod_{t=1}^s oe_{i_t j_t}(Y h_t(X,Y)),
\end{eqnarray*}
where either $i_t=1$ or $j_t=1$ and $h_t(X,Y) \in R[X,Y]$, when
$i_t=1$; $h_t(X,Y) \in I[X,Y]$ when $j_t=1$.
\end{lem}

Proof: We prove the result using induction on $r$, where $\varepsilon$ is
product of $r$ many elementary generators. Let $r=1$ and $\varepsilon
= oe_{pq}(a)$. Note that $a \in R$ when $p=1$, and $a \in I$ when
$q=1$. Given that $oe_{ij}(Xf(X))$ is an elementary generator of
$\EO_{2n}^1(R[X],I[X])$. First we assume $i=1$, hence
$f(X) \in R[X]$.

{\it Case $($1$)$:} Let $(p,q)=(1,j)$. In this case
\begin{eqnarray*}
oe_{1j}(a) ~ oe_{1j}(Y^4X f(Y^4X)) ~ oe_{1j}(-a) &=& oe_{1j}(Y^4X f(Y^4X)).
\end{eqnarray*}

{\it Case $($2$)$:} Let $(p,q)=(1, \sigma(j))$. In this case
\begin{eqnarray*}
oe_{1 \sigma(j)}(a) ~ oe_{1j}(Y^4X f(Y^4X)) ~ oe_{1 \sigma(j)}(-a) &=& oe_{1j}(Y^4X f(Y^4X)).
\end{eqnarray*}

{\it Case $($3$)$:} Let $(p,q)=(1, k), k \ne j , \sigma(j)$. In this case
\begin{eqnarray*}
oe_{1 k}(a) ~ oe_{1j}(Y^4X f(Y^4X)) ~ oe_{1 k}(-a) &=& oe_{1j}(Y^4X f(Y^4X)).
\end{eqnarray*}

{\it Case $($4$)$:} Let $(p,q)=(j, 1))$. In this case
\begin{eqnarray*}
&& oe_{j 1}(a) ~ oe_{1j}(Y^4X f(Y^4X)) ~ oe_{j 1}(-a) \\
&=&  ^{oe_{j1}(a)}  [oe_{1k}(Y^2X f(Y^4X)), oe_{kj}(Y^2)] \\
&=& [oe_{jk}(aY^2X f(Y^4X)) ~ oe_{1k}(Y^2X f(Y^4X)), oe_{k1}(a Y^2) ~ oe_{kj}(Y^2)] \\
&=& oe_{jk}(a Y^2X f(Y^4X)) ~ oe_{1k}(Y^2X f(Y^4X)) ~ oe_{k1}(a Y^2) ~ oe_{kj}(Y^2) \\
&& oe_{1k}(-Y^2X f(Y^4X)) ~ oe_{jk}(- aY^2X f(Y^4X)) ~ oe_{kj}(-Y^2) ~ oe_{k1}(-a Y^2) \\
&=& oe_{jk}(aY^2X f(Y^4X)) ~ oe_{1k}(Y^2X f(Y^4X)) ~ oe_{k1}(a Y^2) ~ oe_{1k}(-Y^2X f(Y^4X)) \\
&& oe_{1k}(Y^2X f(Y^4X)) ~ oe_{kj}(Y^2) ~ oe_{1k}(-Y^2X f(Y^4X)) ~ oe_{kj}(-Y^2) ~ oe_{kj}(Y^2) \\ 
&& [oe_{j1}(-aY), oe_{1k}(Y X f(Y^4X))] ~ oe_{kj}(-Y^2) ~ oe_{k1}(-a Y^2) \\
&=& [oe_{j1}(aY), oe_{1k}(YX f(Y^4X))] ~ oe_{1k}(Y^2X f(Y^4X)) ~ oe_{k1}(aY^2) \\
&& oe_{1k}(Y^2X f(Y^4X))  oe_{1j}(Y^2X f(Y^4X)) [oe_{k1}(-aY^4) oe_{j1}(-aY), \\
&& oe_{1k}(Y X f(Y^4X)) ~ oe_{1j}(Y^3 X f(Y^4X))]  oe_{k1}(-a Y^2) 
\end{eqnarray*}

{\it Case $($5$)$:} Let $(p,q)=(\sigma(j), 1))$. In this case
\begin{eqnarray*}
&& oe_{\sigma(j) 1}(a) ~ oe_{1j}(Y^4X f(Y^4X)) ~ oe_{\sigma(j) 1}(-a) \\
&=& ^{oe_{\sigma(j) 1}(a)} [oe_{1k}(Y^2), oe_{kj}(Y^2X f(Y^4X))] \\
&=& [oe_{\sigma(j) k} (aY^2) ~ oe_{1k}(Y^2), oe_{kj}(Y^2X f(Y^4X))] \\
&=& oe_{\sigma(j) k} (aY^2) ~ oe_{1k}(Y^2) ~ oe_{kj}(Y^2X f(Y^4X)) ~ oe_{1k}(-Y^2) ~ oe_{\sigma(j) k} (-aY^2) \\
&& oe_{kj}(-Y^2X f(Y^4X)) \\
&=& oe_{\sigma(j) k} (aY^2) ~ oe_{1j} (Y^4X f(Y^4X)) ~ oe_{kj}(Y^2X f(Y^4X)) \\
&& [oe_{\sigma(j) 1}(-aY), oe_{1k}(Y)] ~ oe_{kj}(-Y^2X f(Y^4X)) \\
&=& [oe_{\sigma(j) 1}(aY), oe_{1k}(Y)] ~ oe_{1j}  (Y^4X f(Y^4X))\\
&& [oe_{\sigma(j) 1}(-aY), oe_{1k}(Y) oe_{1j}(Y^3X f(Y^4X))]
\end{eqnarray*}

{\it Case $($6$)$:} Let $(p,q)=(k, 1), k \ne j , \sigma(j)$. In this case
\begin{eqnarray*}
&& oe_{k 1}(a) ~ oe_{1j}(Y^4X f(Y^4X)) ~ oe_{k 1}(-a) \\
&=& oe_{kj} (aY^4X f(Y^4X)) ~  oe_{1j}(Y^4X f(Y^4X)) \\
&=& [oe_{k1}(aY^2), oe_{1j}(Y^2X f(Y^4X))] ~ oe_{1j}(Y^4X f(Y^4X)) 
\end{eqnarray*}

Hence the result is true when $i=1$ and $\varepsilon$ is an elementary
generator. Carrying out similar calculation one can show the result is
true when $j=1$ and $\varepsilon$ is an elementary generator. Let us
assume that the result is true when $\varepsilon$ is product of $r-1$
many elementary generators, i.e, $\varepsilon_2 \ldots \varepsilon_r ~
oe_{ij}(Y^{4^{r-1}}X f(Y^{4^{r-1}}(X)) ~ \varepsilon_r^{-1} \ldots
\varepsilon_2^{-1} = \prod_{t=1}^k oe_{p_t q_t}(Y g_t(X,Y))$, where
either $p_t=1$ or $q_t=1$. Note that $g_t(X,Y) \in R[X,Y]$ when
$p_t=1$ and $g_t(X,Y) \in I[X,Y]$ when $q_t=1$.

We now establish the result when $\varepsilon$ is product of $r$ many
elementary generators. We have
\begin{eqnarray*}
&& \varepsilon ~ oe_{ij}(Y^{4^r}X f(Y^{4^r}X)) ~ \varepsilon^{-1}  \\
& = & \varepsilon_1 \varepsilon_2 \ldots \varepsilon_r ~ oe_{ij}(Y^{4^r}X f(Y^{4^r}X)) 
~ \varepsilon_r^{-1} \ldots \varepsilon_2^{-1} \varepsilon_1^{-1} \\
& = & \varepsilon_1 ~ \big( \prod_{t=1}^k oe_{p_t q_t}(Y^4 g'_t(X,Y)) \big) ~ \varepsilon_1^{-1}  \\
& = & \prod_{t=1}^k \varepsilon_1 ~ oe_{p_t q_t}(Y^4 g'_t(X,Y)) ~ \varepsilon_1^{-1} \\
& = & \prod_{t=1}^s oe_{i_t j_t}(Y h_t(X,Y)). 
\end{eqnarray*}

To get the last equality one needs to repeat the calculation which was
done for a single elementary generator. Note that at the last line
either $i_t=1$ or $j_t=1$. Also, note that $h_t(X,Y) \in R[X,Y]$, when
$i_t=1$ and $h_t(X,Y) \in I[X,Y]$, when $j_t=1$ 
\hfill{$\square$}

\begin{nt}
{\rm Let $M$ be a finitely presented $R$-module and $a$ be a
  non-nilpotent element of $R$. Let $R_a$ denote the ring $R$
  localised at the multiplicative set $\{a^i : i \ge 0 \}$ and $M_a$
  denote the $R_a$-module $M$ localised at $\{a^i : i \ge 0 \}$. Let
  $\alpha(X)$ be an element of $\End(M[X])$. The localization map $i:
  M \to M_a$ induces a map $i^*: \End(M[X]) \to \End(M[X]_a)
  = \End(M_a[X])$. We shall denote $i^*(\alpha(X))$ by $\alpha(X)_a$
  in the sequel.}
\end{nt}

We need the following two lemmas.

\begin{lem} \label{equal-auto}
  Let $M$ be a finitely presented $R$-module and $I$ be an ideal of
  $R$. Let $\alpha(X), \beta(X) \in
  {\End}(M[X],IM[X])=ker({\End}(M[X]) \lra {\End}(M[X]/IM[X]))$, with
  $\alpha(0)=\beta(0).$ Let $a$ be a non-nilpotent element in
  $R$ such that $\alpha(X)_a = \beta(X)_a$ in
  ${\End}(M_a[X],IM_a[X])$. Then $\alpha(a^N X)=\beta(a^N X)$ in
  ${\End}(M[X],IM[X])$, for $N \gg 0$.
\end{lem}

\begin{lem} \label{E1-dil-strong} Let $R$ be a commutative ring and 
$I$ be an ideal of $R$. Let $n \ge 3$. Let $a$ be a
  non-nilpotent element in $R$ and $\alpha(X)$ be in
  $\EO_{2n}^1(R_a[X],I_a[X])$, with $\alpha(0)=Id$. Then there exists
  $\alpha^*(X) \in \EO_{2n}^1(R[X],I[X])$, with $\alpha^*(0) = Id.$, such
  that $\alpha^*(X)$ localises to $\alpha(bX)$, for $b \in (a^d)$, $d
  \gg 0$.
\end{lem}

The proofs as in Lemma 3.3 in \cite{cr2} and in Lemma 3.4 in \cite{cr2} work verbatim for Lemma \ref{equal-auto} and Lemma \ref{E1-dil-strong} as above, respectively, and hence the proofs are omitted.

The following result was proved in \cite{acr}. We next apply Lemma \ref{E1-dil-strong} to record a different proof.  

\begin{thm} \label{rel-dil-strong} Let $R$ be a commutative ring and 
$I$ be an ideal of $R$. Let $n \ge 3$. Let $a$ be a
  non-nilpotent element in $R$ and $\alpha(X)$ be in
  $\EO_{2n}(R_a[X],I_a[X])$, with $\alpha(0)=Id$. Then there exists
  $\alpha^*(X) \in \EO_{2n}(R[X],I[X])$, with $\alpha^*(0) = Id.$, such that
  $\alpha^*(X)$ localises to $\alpha(bX)$, for $b \in (a^d)$, $d \gg
  0$.
\end{thm}

Proof: Follows from Lemma \ref{vanderk1} and Lemma \ref{E1-dil-strong}.
\hfill{$\square$}

\section{\large{Orthogonal Modules and Orthogonal Transvections}}

In this section we prove Theorem \ref{eql-rel} which is the main result of this paper. We begin with a sequence of definitions.

\begin{de}
{\rm
Let $M$ be an $R$-module. A {\it bilinear form} on $M$ is a function 
$\beta: M \times M \lra R$ such that $\beta(x, y)$ is $R$-linear as a function
of $x$ for fixed $y$, and $R$-linear as a function of $y$ for fixed $x$. The pair
$(M, \beta)$ is called bilinear form module over $R$. $\beta$ is called an {\it inner product}
if it satisfies non-degeneracy condition, i.e, the natural map induced by $\beta$
 from $P \lra P^*$ is an isomorphism. In this case
the pair $(M, \beta)$ is called inner product module over $R$. A bilinear form
or inner product $\beta$ is called {\it symmetric} if $\beta(x, y) = \beta(y, x)$, for all $x, y \in M$. An inner
product module $(M, \beta)$ will be called {\it inner product space} if $M$ is finitely
generated and projective over $R$.
}
\end{de}

\begin{de} 
{\rm 
An {\it orthogonal $R$-module} is a pair $(P,\langle , \rangle)$,
where $P$ is a finitely generated projective $R$-module of even rank
and $\langle , \rangle: P \times P \lra R$ is a non-degenerate symmetric bilinear form.
This is also known as {\it symmetric inner product space}.}
\end{de}

\begin{de}
{\rm
Let $(P_1,\langle , \rangle_1)$ and $(P_2,\langle , \rangle_2)$ be
two  orthogonal $R$-modules. Their {\it orthogonal sum} is the pair
$(P,\langle , \rangle)$, where $P=P_1 \oplus P_2$ and the inner product
is defined by $\langle (p_1,p_2),(q_1,q_2)\rangle = \langle
p_1,q_1\rangle_1 + \langle p_2,q_2\rangle_2$.}
\end{de}

\begin{de}
{\rm There is a non-degenerate symmetric bilinear form $\langle , \rangle$ on the
$R$-module $R \oplus R^*$, namely $\langle (a_1,f_1),
(a_2,f_2) \rangle = f_2(a_1) + f_1(a_2)$. The orthogonal module $R \oplus R^*$ with this
symmetric bilinear form is denoted by $\mathbb{H}(R)$ and called {\it hyperbolic plane}.
Note that $\mathbb{H}^n(R)$ is the orthogonal sum of $n$-copies of $\mathbb{H}(R)$.}
\end{de}

\begin{de}
{\rm An {\it isometry} of an orthogonal module $(P,\langle , \rangle)$ is
  an automorphism of $P$ which fixes the bilinear form. The group of
  isometries of $(P, \langle , \rangle)$ is denoted by
  $\O(P)$.  }
\end{de}

\begin{de} 
{\rm 
Let $(P, \langle, \rangle)$ be an orthogonal module. H. Bass defined 
{\it orthogonal transvection}
of an orthogonal module $(P, \langle, \rangle)$ is an automorphism of the form
\begin{eqnarray*}
\tau(p) &=& p - \langle u , p \rangle v + \langle v , p \rangle u,
\end{eqnarray*}
where $u,v \in P$ are isotropic, i.e, $\langle u, u \rangle = \langle v, v \rangle = 0$
with $\langle
u,v \rangle=0$, and either $u$ or $v$ is unimodular. It is easy to
check that $\langle \tau(p), \tau(q) \rangle = \langle p, q
\rangle$, i.e, $\tau \in \O(P)$ and $\tau$ has an inverse $\sigma(p) = p + \langle u, p
\rangle v - \langle v, p \rangle u$.

The subgroup of $\O(P)$ generated by the orthogonal
transvections is called orthogonal transvection group and denoted by $\Trans_{\O}(P, \langle , \rangle)$ (see
\cite{bass2} or \cite{HO}).  }
\end{de}

{\bf Now onwards $Q$ will denote $(R^2 \oplus P)$ with induced form
  on $(\mathbb{H}(R)~\oplus~P)$, and $Q[X]$ will denote $(R[X]^2 \oplus
  P[X])$ with induced form on $(\mathbb{H}(R[X])~\oplus~P[X])$.}

\begin{de} {\rm The orthogonal transvections of $Q=(R^2 \oplus P)$ of 
the form
    \begin{eqnarray*}
      (a, b, p) & \mapsto & (a, b + \langle p, q \rangle, p-aq),
     \end{eqnarray*} 
or of the form 
    \begin{eqnarray*}
      (a, b, p) & \mapsto & (a + \langle p, q \rangle, b,  p-bq), 
    \end{eqnarray*} 
where $a, b \in R$ and $p, q \in P$, 
are called {\it elementary orthogonal transvections}. Let us denote the first 
isometry by $\rho(q)$ and the second one by $\mu(q)$. It can 
be verified that the elementary orthogonal transvections are orthogonal 
transvections on $Q$. Indeed, consider $(u, v) =((0,1,0),(0,0,q))$  to get 
$\rho(q)$ and consider $(u, v)= ((1,0,0),(0,0,q))$ to get $\mu(q, \beta)$.

The subgroup of $\Trans_{\O}(Q, \langle , \rangle)$ generated by
elementary orthogonal transvections is denoted by $\ETrans_{\O}(Q,
\langle , \rangle)$.}
\end{de}

\begin{de} 
{\rm Let $I$ be an ideal of $R$. The group of {\it relative orthogonal
    transvections} to an ideal  $I$ is generated by the orthogonal
  transvections of the form $\sigma(p) = p - \langle u, p \rangle v +
  \langle v, p \rangle u$, where either $u \in IP$ or $v \in IP$. The group 
  generated by relative orthogonal
  transvections is denoted by $\Trans_{\O}(P,IP, \langle , \rangle)$.}
\end{de}

\begin{de} {\rm Let $I$ be an ideal of $R$. The elementary orthogonal
transvections of $Q$ of the form $\rho(q), \mu(q)$, where
$q \in IP$ are called {\it relative elementary 
orthogonal transvections} to an ideal $I$.

The subgroup of $\ETrans_{\O}(Q, \langle , \rangle)$ generated by
relative elementary orthogonal transvections is denoted by
$\ETrans_{\O}(IQ,\langle , \rangle )$. The normal closure of
$\ETrans_{\O}(IQ,\langle , \rangle )$ in
$\ETrans_{\O}(Q, \langle , \rangle)$ is denoted by
$\ETrans_{\O}(Q,IQ, \langle , \rangle)$.}
\end{de}

\begin{rmk} \label{free} 
Let $P=\oplus_{i=1}^{2n} Re_i$ be a free $R$-module with $R=2R$. 
The non-degenerate symmetric bilinear form $\langle,\rangle$ on $P$
corresponds to a symmetric matrix $\varphi$ with
respect to the basis $\{e_1, e_2, \ldots, e_{2n} \}$ of $P$ and we write
$\langle p, q \rangle = p^t \varphi q$.

  In this case the orthogonal transvection $\tau(p) = p - \langle u,
  p \rangle v + \langle v, p \rangle u$ corresponds to the matrix 
  $(I_{2n} - v u^t \varphi  + u v^t \varphi)$ and the group generated by them
  is denoted by $\Trans_{\O}(P, \langle , \rangle_{\varphi})$.

Also in this case $\ETrans_{\O}(Q, \langle , \rangle_{\widetilde{\psi}_1 \perp \varphi})$ will
be generated by the matrices of the form $\rho_{\varphi}(q) =
\Big( \begin{smallmatrix} 1 & 0 & 0 \\ 0 & 1 & q^t \varphi \\ -q &
  0 & I_{2n} \end{smallmatrix} \Big)$, and $\mu_{\varphi}(q) =
\Big( \begin{smallmatrix} 1 & 0 & q^t \varphi \\ 0 & 1 & 0 \\ 0 &
-q & I_{2n} \end{smallmatrix} \Big)$.

Note that for standard symmetric matrix $\widetilde{\psi}_n$ and for $q=(q_1, \ldots, q_{2n}) \in 
R^{2n}$ with $q^t \widetilde{\psi}_n q = 0$, we have
\begin{eqnarray} \label{relatn5}
\rho_{\widetilde{\psi}_n}(q) &=& \prod_{i=3}^{2n+2} oe_{i1}(-q_{i-2}), \\ 
\label{relatn6}
\mu_{\widetilde{\psi}_n}(q) &=& \prod_{i=3}^{2n+2} oe_{1i}(-q_{\sigma(i-2)}).  
\end{eqnarray}
\end{rmk}

In the following four lemmas we shall use the assumptions and notations in the
  statement of Remark \ref{free}.

\begin{lem} \label{kopeiko} 
Let $R$ be a commutative ring with $R=2R$, and $I$ be an ideal of $R$. Let $(P, \langle, \rangle)$
be an orthogonal $R$-module with $P$ free $R$-module of rank $2n$, $n \ge 2$ and $Q= R^2 
\oplus P$ with the induced form on $\mathbb{H}(R) \oplus P$. If the symmetric bilinear form $\langle, \rangle$ correspond 
(w.r.t. some basis) to  $\widetilde{\psi}_n$, the standard symmetric matrix, then $\Trans_{\O}(Q,IQ, 
\langle , \rangle_{\widetilde{\psi}_{n+1}}) = \EO_{2n+2}(R,I)$.
\end{lem}

Proof: For proof see \S 2 of \cite{SK}.

\begin{lem} \label{free,psi} Let $R$ be a commutative ring with
  $R=2R$, and let $I$ be an ideal of $R$. Let $(P, \langle, \rangle)$
be an orthogonal $R$-module with $P$ free $R$-module of rank $2n$, $n \ge 2$ and $Q= R^2 
\oplus P$ with the induced form on $\mathbb{H}(R) \oplus P$. If the symmetric bilinear form $\langle, \rangle$ correspond 
(w.r.t. some basis) to  $\widetilde{\psi}_n$, the standard
  symmetric matrix, then $\ETrans_{\O}(Q,IQ,
  \langle , \rangle_{\widetilde{\psi}_{n+1}}) = \EO_{2n+2}(R,I)$.
\end{lem}

Proof: We first show $\ETrans_{\O}(Q,IQ, \langle ,
\rangle_{\widetilde{\psi}_{n+1}})$ is a subset of $\EO_{2n+2}(R,I)$. An element of $\ETrans_{\O}(Q,IQ, \langle ,
\rangle_{\widetilde{\psi}_{n+1}})$ is of the form $T_1(q) T_2(s) T_1(q)^{-1}$, where $q \in R^{2n}$, $s \in I^{2n} (\subseteq R^{2n})$. Here $T_1$ and $T_2$ can be either of $\rho_{\widetilde{\psi}_n}$ or $\mu_{\widetilde{\psi}_n}$.
Using equations (\ref{relatn5}) and (\ref{relatn6}) we
show that either of the above elements belong to $\EO_{2n+2}(R,I)$,
 and hence $\ETrans_{\O}(Q,IQ, \langle , \rangle_{\widetilde{\psi}_{n+1}}) \subseteq
\EO_{2n+2}(R,I)$.

To show the other inclusion we recall that $\EO_{2n+2}(R,I)$ is generated
by the elements $g ~ oe_{ij}(x) g^{-1}$, where $g \in \EO_{2n+2}(R), x \in I$, and
either $i=1$ or $j=1$ (see Lemma \ref{equiv-defn}). Using commutator
relation $[oe_{ik}(a), oe_{kj}(b)] = oe_{ij}(ab)$ and the 
equations (\ref{relatn5}), (\ref{relatn6}) we can show that $\EO_{2n+2}(R,I) 
\subseteq \ETrans_{\O}(Q,IQ, \langle , \rangle_{\widetilde{\psi}_{n+1}})$, and hence
the equality is established.
\hfill{$\square$}

\begin{lem} \label{phi=phi*,ab}
Let $P$ be a free $R$-module of rank $2n$. Let $(P,\langle ,
\rangle_{\varphi})$ and $(P,\langle , \rangle_{\varphi^*})$ be two
orthogonal $R$-modules with $\varphi= \varepsilon^t 
\varphi^*  \varepsilon$, for some $\varepsilon \in
\GL_{2n}(R)$. Then
\begin{eqnarray*} 
\Trans_{\O}(P,\langle , \rangle_{\varphi}) &=& \varepsilon^{-1} ~ 
\Trans_{\O}(P,\langle , \rangle_{\varphi^*}) ~ \varepsilon,\\ 
\ETrans_{\O}(Q, \langle , \rangle_{\widetilde{\psi}_1 \perp \varphi}) &=& (I_2 \perp
\varepsilon)^{-1}~ \ETrans_{\O}(Q,\langle , \rangle_{\widetilde{\psi}_1 \perp \varphi^*}) ~ (I_2
\perp \varepsilon).
\end{eqnarray*}
\end{lem}

Proof: In the free case for orthogonal transvections we have
\begin{eqnarray*}
(I_{2n} - v u^t \varphi + u v^t \varphi) & = & \varepsilon^{-1} ~ (I_{2n} - \tilde{v} \tilde{u}^t \varphi^*
 - \tilde{u} \tilde{v}^t \varphi^*) ~ \varepsilon,
\end{eqnarray*}
where $\tilde{u}=\varepsilon u$ and $\tilde{v} = 
\varepsilon v$. Hence the first equality follows. 

For elementary orthogonal transvections we have 
\begin{eqnarray*}
(I_2 \perp \varepsilon)^{-1} \rho_{\varphi^*}(q) (I_2 \perp \varepsilon) &=& \rho_{\varphi} (\varepsilon^{-1} q),\\ 
(I_2 \perp \varepsilon)^{-1} \mu_{\varphi^*}(q) (I_2 \perp \varepsilon) &=& \mu_{\varphi} (\varepsilon^{-1} q),
\end{eqnarray*}
hence the second equality follows.
\hfill{$\square$}

\begin{lem} \label{phi=phi*,rel} Let $I$ be an ideal of $R$ and $P$ be
  a free $R$-module of rank $2n$. Let $(P,\langle , \rangle_{\varphi})$
  and $(P,\langle , \rangle_{\varphi^*})$ be two orthogonal $R$-modules
  with $\varphi= \varepsilon^t  \varphi^* \varepsilon$, for some $\varepsilon \in \GL_{2n}(R)$. Then
\begin{eqnarray*} 
\Trans_{\O}(P,IP,\langle , \rangle_{\varphi}) &=& 
\varepsilon^{-1} ~ \Trans_{\O}(P,IP,\langle , \rangle_{\varphi^*}) ~ \varepsilon,\\ 
\ETrans_{\O}(Q,IQ,\langle , \rangle_{\widetilde{\psi}_1 \perp \varphi}) &=& (I_2 \perp
\varepsilon)^{-1} ~ \ETrans_{\O}(Q,IQ,\langle , \rangle_{\widetilde{\psi}_1 \perp \varphi^*}) ~ (I_2
\perp \varepsilon).
\end{eqnarray*}
\end{lem}

Proof: Using the three equations appearing in the proof of Lemma
\ref{phi=phi*,ab}, we get these equalities.
\hfill{$\square$}

\begin{pr} \label{local-case,rel} Let $R$ be a commutative ring with
  $R=2R$, and let $I$ be an ideal of $R$. Let $(P,\langle ,
  \rangle_{\varphi})$ be an orthogonal $R$-module with $P$ free of rank
  $2n$, $n \ge 2$ and $Q = R^2 \oplus P$ with the induced form on $\mathbb{H}(R) \oplus P$. If $\varphi=\varepsilon^t  \widetilde{\psi}_n  \varepsilon$, 
  for some $\varepsilon \in \GL_{2n}(R)$, then
  $\Trans_{\O}(Q,IQ,\langle , \rangle_{\widetilde{\psi}_1 \perp \varphi}) =
  \ETrans_{\O}(Q,IQ,\langle , \rangle_{\widetilde{\psi}_1 \perp \varphi})$.
\end{pr}

Proof: Using Lemma \ref{kopeiko}, Lemma \ref{free,psi}, and Lemma
\ref{phi=phi*,rel} we get,

\begin{eqnarray*}
\Trans_{\O}(Q,IQ,\langle , \rangle_{\widetilde{\psi}_1 \perp \varphi}) &=& (I_2
\perp \varepsilon)^{-1}~ \Trans_{\O}(Q,IQ,\langle , \rangle_{\widetilde{\psi}_{n+1}})
  ~ (I_2 \perp \varepsilon)\\ 
  &=& (I_2 \perp \varepsilon)^{-1}~
  \EO_{2+2n}(R,I) ~ (I_2 \perp \varepsilon),
\end{eqnarray*}
and
\begin{eqnarray*}
\ETrans_{\O}(Q,IQ,\langle , \rangle_{\widetilde{\psi}_1 \perp \varphi}) &=& (I_2 \perp
\varepsilon)^{-1} ~ \ETrans_{\O}(Q,IQ,\langle , \rangle_{\widetilde{\psi}_{n+1}}) ~ (I_2
\perp \varepsilon)\\ 
&=& (I_2 \perp \varepsilon)^{-1} ~ \EO_{2+2n}(R,I)
~ (I_2 \perp \varepsilon),
\end{eqnarray*}
and hence the equality is established. 
\hfill{$\square$}

\begin{de}
\rm{
An orthogonal module $(P, \langle, \rangle)$ over the ring $R$ is called {\it split} if there
exists a submodule $N \subseteq P$ such that $N$ is a direct summand of $P$
and such that $N$ is precisely equal to its orthogonal complement $N^{\perp} = \{ p \in P :
\langle p, n \rangle = 0 ~{\rm for~ all}~ n \in N \}$. 

Moreover, an orthogonal module $(P, \langle, \rangle)$ over the ring $R$ is called {\it locally split} if 
$(P_\gm, \langle, \rangle)$ 
is a split orthogonal $R_\gm$-module for every maximal ideal $\gm$ of $R$.
}
\end{de}

\begin{lem} \label{split&free} (See Lemma 6.3, Chapter I in \cite{MH})
Let $R$ be a ring such that  every finitely generated projective module over $R$ is free. Then 
an inner product space over $R$ is split if and only if it possesses a basis so that the associated 
inner product matrix has the form $\big( \begin{smallmatrix} 0 & I \\ I & A \end{smallmatrix} \big)$.
If we also assume that $2$ is a unit in the ring, then every split inner product space has matrix
$\big( \begin{smallmatrix} 0 & I \\ I & 0 \end{smallmatrix} \big)$ with respect to a suitable basis.
\end{lem}

\begin{rmk} \label{localcase} In view of Proposition \ref{local-case,rel} 
and above lemma, for any
    split orthogonal module $(P,\langle , \rangle_{\varphi})$ over a local
    ring $(R,\gm)$ with $R=2R$, we have $\Trans_{\O}(Q,IQ,\langle,\rangle_{\widetilde{\psi}_1
      \perp \varphi} ) = \ETrans_{\O}(Q,IQ, \langle , \rangle_{\widetilde{\psi}_1
      \perp \varphi} )$. Here $I$ is an ideal of the ring $R$.
\end{rmk}

Next we establish dilation principle for relative elementary orthogonal
transvection group.

\begin{lem} \label{rel-dilatn-ETrans_O} Let $R$ be a commutative ring
  with $R=2R$, and let $I$ be an ideal of $R$. Let $(P,\langle ,
  \rangle)$ be an orthogonal $R$-module with rank of $P$ is
   $2n$, $n \ge 2$, and $Q=R^2 \oplus P$ with the induced form on $\mathbb{H}(R) \oplus P$. Suppose that $a$ is a
  non-nilpotent element of $R$ such that $P_a$ is a free $R_a$ module,
  $(P_a, \langle, \rangle)$ is split orthogonal $R_a$-module, 
  and the bilinear form $\langle, \rangle$ corresponds to the
  symmetric matrix $\varphi$ (w.r.t. some basis). Let
  $\alpha(X) \in \Aut(Q[X])$, with $\alpha(0) = Id$, and $\alpha(X)_a 
  \in \ETrans_{\O}(Q_a[X],IQ_a[X],\langle , \rangle_{\widetilde{\psi}_1 \perp
    \varphi})$. Then, there exists $\alpha^*(X) \in 
    \ETrans_{\O}(Q[X],IQ[X],\langle , \rangle)$, with $\alpha^*(0) =
  Id.$, such that $\alpha^*(X)$ localises to $\alpha(bX)$, for $b \in
  (a^d)$, $d \gg 0$.
\end{lem}

Proof: We have $P_a \cong R_a^{2n}$. Let $e_1, \ldots, e_{2n+2}$ be the standard basis of $Q_a$ with respect to which the bilinear form on $Q_a$ will correspond to $\widetilde{\psi}_1 \perp \varphi$. Since $(P_a, \langle, \rangle)$ is a split orthogonal $R$-module with $R_a = 2R_a$, we have $\varphi = \varepsilon^t \widetilde{\psi}_n \varepsilon$, for some $\varepsilon \in \GL_{2n}(R_a)$ by Lemma \ref{split&free}. Therefore, 
$\ETrans_{\O}(Q_a[X], IQ_a[X], \langle, \rangle_{\widetilde{\psi}_1 \perp \varphi}) = (I_2 \perp \varepsilon)^{-1} \EO_{2n+2}(R_a[X], I_a[X]) (I_2 \perp \varepsilon)$ by Lemma \ref{free,psi}, and Lemma
\ref{phi=phi*,rel}. Hence, $\alpha(X)_a = (I_2 \perp \varepsilon)^{-1} ~\beta(X)~ (I_2 \perp \varepsilon)$, for some $\beta(X) \in \EO_{2n+2}(R_a[X], I_a[X])$, with $\beta(0)=Id.$ By Lemma \ref{vanderk1} we have
\begin{eqnarray*}
\EO_{2n+2}(R_a[X], I_a[X]) &=& \EO_{2n+2}^1(R_a[X], I_a[X]) ~\cap  \O_{2n+2}(R_a[X], I_a[X]).
\end{eqnarray*} 
Hence we can write $\beta(X) = \prod_t \gamma_t ~ oe_{i_t j_t}(X f_t(X)) ~ \gamma_t^{-1}$, where either $i_t =1$, or $j_t=1$, and $\gamma_t \in \EO_{2n+2}^1(R_a, I_a)$. Note that $f_t(X) \in R_a[X]$, when $i_t=1$ and $f_t(X) \in I_a[X]$, when $j_t=1$. Using Lemma \ref{ness4-E1} we get $\beta(Y^{4^r}X) = \prod_k oe_{i_k j_k}(Y h_k(X,Y)/ a^m)$, with either $i_k =1$ or $j_k=1$. Note that $h_k(X,Y) \in R[X,Y]$, when $i_k=1$ and $h_k(X,Y) \in I[X,Y]$, when $j_k=1$. We have 
\begin{eqnarray*}
oe_{1 j_k}(Y h_k(X,Y)/ a^m) &=& I_{2n+2} - (-Y h_k(X,Y)/ a^m) ~ e_1 ~ e_{\sigma(j_k)}^t  ~ \widetilde{\psi}_{n+1} \\
&& + (-Y h_k(X,Y)/ a^m) ~ e_{\sigma(j_k)} ~ e_1^t ~ \widetilde{\psi}_{n+1}, ~ {\rm for} ~ j_k \ge 3, \\
oe_{i_k 1}(Y h_k(X,Y)/ a^m) &=& I_{2n+2} - (-Y h_k(X,Y)/ a^m) ~ e_{i_k} ~ e_2^t ~ \widetilde{\psi}_{n+1}\\
&& + (- Y h_k(X,Y)/ a^m) ~ e_2 ~ e_ {i_k}^t ~ \widetilde{\psi}_{n+1}, ~ {\rm for} ~ i_k \ge 3. \\
\end{eqnarray*}

Let $\varepsilon_1, \ldots, \varepsilon_{2n}$ be the columns of the matrix $\varepsilon \in \GL_{2n}(R_a)$. Let $\widetilde{e_i}^t$ denote the column vector $(I_2 \perp \varepsilon) e_i$ of length $2n+2$. Note that $\widetilde{e_1}=e_1, \widetilde{e_2}=e_2$, and $\widetilde{e_i}^t = (0, 0, \varepsilon_{i-2}^t)$, for $i \ge 3$. Using Lemma \ref{phi=phi*,ab}  we can write $\alpha(Y^{4^r}X)_a$ as product of elements of the form
\begin{eqnarray*}
& I_{2n+2} - (-Y h_k(X,Y)/ a^m) \widetilde{e_1} \widetilde{e}_{\sigma(j_k)}^t  \left( \begin{smallmatrix} \widetilde{\psi}_1 & 0 \\ 0 & \varphi \end{smallmatrix} \right ) + (-Y h_k(X,Y)/ a^m)  \widetilde{e}_{\sigma(j_k)} \widetilde{e_1}^t  \left ( \begin{smallmatrix} \widetilde{\psi}_1 & 0 \\ 0 & \varphi \end{smallmatrix} \right )& \\
&= \mu_{\varphi}((Y h_k(X,Y)/ a^m) \varepsilon_{\sigma(j_k) -2}), & \\
& I_{2n+2} - (-Y h_k(X,Y)/ a^m) \widetilde{e}_{i_k} \widetilde{e_2}^t  \left ( \begin{smallmatrix} \widetilde{\psi}_1 & 0 \\ 0 & \varphi \end{smallmatrix} \right ) + (-Y h_k(X,Y)/ a^m) \widetilde{e_2} \widetilde{e}_{i_k}^t  \left ( \begin{smallmatrix} \widetilde{\psi}_1 & 0 \\ 0 & \varphi \end{smallmatrix} \right ) & \\
&= \rho_{\varphi}(-(Y h_k(X,Y)/ a^m) \varepsilon_{i_k -2}), & 
\end{eqnarray*}
 for $i_k, j_k \ge 3$. Note that $\alpha(Y^{4^r}X)_a \in \ETrans_{\O}(Q_a[X, Y],IQ_a[X, Y],\langle , \rangle_{\widetilde{\psi}_1 \perp
    \varphi})$, hence $\alpha(Y^{4^r}X)_a = id ~{\rm mod}~ (IQ_a[X, Y])$. Since $\rho_{\varphi}$ and $\mu_{\varphi}$ satisfy the splitting property $\rho_{\varphi}(q_1 + q_2) = \rho_{\varphi}(q_1) \rho_{\varphi}(q_2)$ and $\mu_{\varphi}(q_1 + q_2) = \mu_{\varphi}(q_1) \mu_{\varphi}(q_2)$, we get $\alpha(Y^{4^r}X)_a$ is product of elements of the form $T_1((Y f_k(X, Y)/a^m) \varepsilon_{p_k})  T_2((Y g_k(X, Y)/a^m) \varepsilon_{q_k}) \\T_1(-(Y f_k(X, Y)/a^m) \varepsilon_{p_k}) $, where $T_1, T_2$ are either $\rho_{\varphi}$ or $\mu_{\varphi}$, $f_k(X, Y) \in R[X, Y]$, $g_k(X, Y) \in I[X, Y]$, and $p_k , q_k \ge 3$.

 Let $s \ge 0$ be an integer such that $\widetilde{\varepsilon_i} = a^s \varepsilon_i \in P$ for all $i = 1, \ldots, 2n$. Let $d = s+m$. Therefore $\alpha((a^{d}Y)^{4^r} X)_a$ is product of elements of the form $T_1((a^d Y f_k(X, a^dY)/a^m) \varepsilon_{p_k})  \\T_2((a^d Y g_k(X, a^d Y)/a^m) \varepsilon_{q_k}) T_1(-(a^d Y f_k(X, a^d Y)/a^m) \varepsilon_{p_k}) $
  Substituting $Y=1$ we get $\alpha(a^dX)_a$ is product of elements of the forms 
\begin{eqnarray*}
T_1(a^s f'_k(X) \varepsilon_{p_k})  T_2(a^s  g'_k(X) \varepsilon_{q_k}) T_1(-(a^m f'_k(X) \varepsilon_{p_k}).
\end{eqnarray*}

Let us set $\alpha^*(X)$ to be the product of elements of the forms 
\begin{eqnarray*}
T_1(f'_k(X) \widetilde{\varepsilon}_{p_k})  T_2(g'_k(X) \widetilde{\varepsilon}_{q_k}) T_1(- f'_k(X) \widetilde{\varepsilon}_{p_k}).
\end{eqnarray*}
From the construction it is clear that $\alpha^*(X)$ belongs to $\ETrans_{\O}(Q[X], IQ[X], \langle, \rangle)$, $\alpha^*(0)=Id.$, and  $\alpha^*(X)$ localises to $\alpha(bX)$, for some $b \in (a^d), d \gg 0$.
\hfill{$\square$}

\begin{lem} \label{LG-ETrans-rel} Let $R$ be a commutative ring with
  $R=2R$, and let $I$ be an ideal of $R$. Let $(P,\langle , \rangle)$
  be a locally split orthogonal $R$-module with $P$ is of rank $2n$, $n \ge 2$,
  and $Q = R^2 \oplus P$ with the induced form on $\mathbb{H}(R) \oplus P$. 
  Let $\alpha(X) \in {\O}(Q[X])$, with $\alpha(0) = Id$.  If $\alpha(X)_\gm \in
  \ETrans_{\O}(Q_\gm[X],IQ_\gm[X],\langle , \rangle_{\widetilde{\psi}_1 \perp
    \varphi_\gm})$, for each maximal ideal $\gm$ of $R$, then 
    $\alpha(X) \in \ETrans_{\O}(Q[X],IQ[X],\langle , \rangle)$.
\end{lem}

Proof: One can suitably choose an element $a_\gm$ from $R \setminus
\gm$ such that $\alpha(X)_{a_\gm}$ belongs to
$\ETrans_{\O}(Q_{a_\gm}[X],IQ_{a_\gm}[X])$. Let us set $\gamma(X,Y) =
\alpha(X+Y) \alpha(Y)^{-1}$. Note that $\gamma(X,Y)_{a_\gm}$ belongs
to $\ETrans_{\O}(Q_{a_\gm}[X,Y],IQ_{a_\gm}[X,Y])$, and $\gamma(0,Y) =
Id$. From Lemma \ref{rel-dilatn-ETrans_O} it follows that $\gamma(b_\gm
X,Y) \in \ETrans_{\O}(Q[X,Y],IQ[X,Y])$, for $b_\gm \in (a_\gm^d)$, where $d
\gg 0$. Note that the ideal generated by $a_\gm^d$'s is the whole ring
$R$. Therefore, $c_1 a_{\gm_1}^d+ \ldots + c_k a_{\gm_k}^d = 1 $, for
some $c_i \in R$. Let $b_{m_i}= c_i a_{m_i}^d \in (a_{m_i}^d)$. It is
easy to see that $\alpha(X)=\prod_{i=1}^{k-1}\gamma(b_{m_i}X,T_i)
\gamma(b_{m_k}X,0)$, where $T_i = b_{m_{i+1}}X+ \cdots +
b_{m_k}X$. Each term in the right hand side of this expression belongs
to $\ETrans_{\O}(Q[X], IQ[X])$ and hence $\alpha(X) \in \ETrans_{\O}(Q[X],
IQ[X])$.
\hfill{$\square$}

\medskip
We now establish equality of the orthogonal transvection group and the
elementary orthogonal transvection group (in the relative case to an 
ideal) of a locally split orthogonal $R$-module
with $R=2R$. An absolute version of this result (i.e, when $I=R$)
was proved in \cite{bbr} (see Theorem 3.10). Before proving the main
result we establish a lemma
to show that orthogonal transvections are homotopic to identity.

\begin{lem} \label{ortho-trans-hom-to-identity}
Let $(P, \langle, \rangle)$ be an orthogonal $R$-module and $\alpha \in
\Trans_{\O}(P, \langle, \rangle)$. Then there exists $\beta(X) \in
\Trans_{\O}(P[X], \langle, \rangle)$ such that $\beta(1)=\alpha$ and
$\beta(0)=Id.$
\end{lem}

Proof: As $\alpha \in \Trans_{\O}(P, \langle, \rangle)$, it is
product of orthogonal transvections of the form $\tau$, where
$\tau$ takes $p \in P$ to $p - \langle u, p \rangle v + \langle v, p
\rangle u$, where $u, v \in P$
are isotropic with $\langle u, v \rangle = 0$, and either $u$ or
$v$ is unimodular. Define $\tau X$ as the map which takes $p \in P$
to either $p - \langle u, p \rangle vX + \langle vX, p \rangle u$ or 
$p - \langle uX, p \rangle v + \langle v, p \rangle uX$. This choice depends on
whether $u$ is unimodular or $v$ is unimodular. Note that $uX$
represents $u$ {\it times} $X$ and $vX$ represents $v$ {\it times} $X$.
 Also, note that $uX, vX \in P[X]$. We set $\beta(X)$ to be the product of
elements of the form $\tau X$, whenever $\tau$ appears in the
expression of $\alpha$. Then $\beta(1)=\alpha$ and $\beta(0)=Id$.
\hfill{$\square$}

\begin{thm} \label{eql-rel} Let $R$ be a commutative ring with $R=2R$,
  and let $I$ be an ideal of $R$. Let $(P, \langle , \rangle)$ be a
 locally split orthogonal $R$-module with $P$ is of rank $2n$, $n \ge 2$,
 and $Q = R^2 \oplus P$ with the induced form on $\mathbb{H}(R) \oplus P$.
 Then $\Trans_{\O}(Q,IQ,\langle , \rangle) = \ETrans_{\O}(Q,IQ,\langle , \rangle)$.
\end{thm}

Proof: We have $\ETrans_{\O}(Q,IQ,\langle,\rangle) \subseteq
\Trans_{\O}(Q,IQ,\langle,\rangle)$. We need to show other inclusion.
Let us choose $\alpha$ from $\Trans_{\O}(Q,IQ,\langle,\rangle)$.  By
Lemma \ref{ortho-trans-hom-to-identity} there exists $\alpha(X)$ in
$\Trans_{\O}(Q[X],IQ[X],\langle,\rangle)$ such that $\alpha(1) =
\alpha$ and $\alpha(0) = Id$. Note that
$\Trans_{\O}(Q_{\gm}[X],IQ_{\gm}[X],\langle,\rangle_{\widetilde{\psi}_1 \perp
  \varphi_{\gm}}) = \ETrans_{\O}(Q_{\gm}[X],IQ_{\gm}[X],\langle,
\rangle_{\widetilde{\psi}_1 \perp \varphi_{\gm}})$, for each maximal ideal $\gm$
of $R$ (follows from Remark \ref{localcase}). Hence $\alpha(X)_\gm$
belongs to
$\ETrans_{\O}(Q_{\gm}[X],IQ_{\gm}[X],\langle,\rangle_{\widetilde{\psi}_1 \perp
  \varphi \otimes R_{\gm}[X]})$, for each maximal ideal $\gm$ of
$R$. Therefore, $\alpha(X) \in
\ETrans_{\O}(Q[X],IQ[X],\langle,\rangle)$ (see Lemma
\ref{LG-ETrans-rel}). Substituting $X=1$ we get the result.
\hfill{$\square$}

\medskip
In closing we make Remark \ref{final-rmk} below
for which  we need the following elementary observation.   

\begin{lem}  \label{local-split}
Let $(P, \langle, \rangle)$ be a split orthogonal $R$-module. Then $(P_\gm, \langle, \rangle)$ 
is a split orthogonal $R_\gm$-module for every maximal ideal $\gm$ of $R$.
\end{lem}

Proof: Let us consider an equivalent form of the definition of split orthogonal $R$-modules
as it is stated in \S 6, Chapter I, \cite{MH}. The orthogonal module $(P, \langle, \rangle)$ is 
split if it is direct sum of two submodules $M$ and $N$ which are dually paired to $R$ by the 
inner product,
\begin{eqnarray*}
M \stackrel{\cong}{\lra} \Hom_R(N,R), ~{\rm and}~ N \stackrel{\cong}{\lra} \Hom_R(M,R)
\end{eqnarray*}
and such that $N$ is self orthogonal, i.e, $\langle N, N \rangle= 0$. Tensoring with $R_\gm$ we get
$P_\gm = M_\gm \oplus N_\gm$.  Moreover, $P$ projective will imply both $M$ and $N$ are 
projective and hence finitely presented ($M$ finitely presented means there exists an exact 
sequence $R^k \lra R^l \lra M$, for suitable natural numbers $k,l$). Therefore, by 
Proposition 2.13" in Chapter I, \cite{lam} we get
\begin{eqnarray*}
M_\gm \stackrel{\cong}{\lra} \Hom_{R_\gm}(N_\gm, R_\gm), ~{\rm and}~ N_\gm  
\stackrel{\cong}{\lra} \Hom_{R_\gm}(M_\gm,R_\gm)
\end{eqnarray*}
Also, $N$ is self orthogonal will imply $N_\gm$ is self-orthogonal and hence 
$(P_\gm, \langle, \rangle)$ is a split orthogonal $R_\gm$-module for every maximal ideal $\gm$ 
of $R$.
\hfill{$\square$}

\begin{rmk}\label{final-rmk} In view of the above lemma the result as in Theorem \ref{eql-rel} holds when  
$(P, \langle , \rangle)$ is assumed to  be a split orthogonal $R$-module.
\end{rmk}

\medskip

\noindent {Acknowledgement:} The author thanks Department of 
Science and Technology, Govt. of India for INSPIRE Faculty Award 
[IFA-13 MA-24] that supported this work.


\begin{thebibliography}{VAST}
\setlength {\itemsep}{-.1ex}


\bibitem{acr} H. Apte, P. Chattopadhyay, R.A. Rao, A Local Global
theorem for extended ideals, J. Ramanujan Math. Soc. 27, No. 1 (2012), 1--20.



\bibitem{bbr} A. Bak, R. Basu, R.A. Rao, Local-Global principle for
  transvection groups, Proceedings of the American Mathematical
  Society 138 no. 4 (2010) 1191--1204






\bibitem{bass2} H. Bass, Unitary algebraic $K$-theory, Lecture Notes
in Mathematics 343 (1973) 57--265


\bibitem{cr} P. Chattopadhyay, R.A. Rao, Elementary symplectic orbits
  and improved $K_1$-stability, Journal of $K$-Theory 7 (2011)
  389--403.
  
\bibitem{cr2}  P. Chattopadhyay, R.A. Rao,  Equality of elementary and
symplectic orbits,  J. Pure Appl. Algebra 219 no. 12 (2015) 5363--5386.
 

\bibitem{HO} A.J. Hahn, O.T. O'Meara, The classical groups and K-Theory,
Grundlehren der Mathematischen Wissenschaften (A Series of Comprehensive Studies in Mathematics),
vol. 291, Springer-Verlag, 1989.


\bibitem{vdK1} W. van der Kallen, A group structure on certain orbit
  sets of unimodular rows, Journal of Algebra 82 no. 2 (1983) 363--397.



\bibitem{lam} T.Y. Lam, Serre's Problem on Projective Modules, Springer
Monographs in Mathematics, Springer, 2006.


\bibitem{MH} J. Milnor, D. Husemoller, Symmetric bilinear forms, Ergebnisse
der Mathematik und ihrer Grenzgebiete (A Series of Modern Surveys in
Mathematics), vol. 73, Springer-Verlag, 1973.


\bibitem{SV} A.A. Suslin, L.N. Vaserstein, Serre's problem on
  Projective Modules over Polynomial Rings and Algebraic K-theory,
  Math.  USSR Izvestija 10 (1976) 937--1001.





\bibitem{SK} A.A. Suslin, V.I. Kope\u\i ko, Quadratic Modules and
  Orthogonal Group over Polynomial Rings (Russian), Algebraic numbers
  and finite groups, Zap. Nau\v cn. Sem. Lenin grad Otdel. Mat. Inst.
  Steklov. (LOMI ), 81 (123) No. 3 1970 328--351.

\bibitem{V3} L.N.  Vaserstein, On the normal subgroups of ${\rm
   GL}_{n}$ over a ring,  Algebraic $K$-theory, Evanston 1980 (Proc.
  Conf., Northwestern Univ., Evanston, Ill., 1980) 456--465; Lecture
  Notes in Math. 854, Springer, Berlin-New York, 1981.
 

\end{thebibliography}
\end{document}